\documentclass[11pt]{amsart}
\usepackage{amsmath,amssymb,enumerate}
\usepackage{epsf,epsfig,amsfonts,graphicx,color}  
 \usepackage[applemac]{inputenc} 

 \date{November, 2010}   
  
 \numberwithin{equation}{section}   
\DeclareMathOperator{\const}{const}%

\newtheorem{Th}{\rm\bf Theorem}

\theoremstyle{definition}

\newtheorem{remark}{\rm\bf Remark}
\newtheorem{example}{\rm\bf Example}
\newcommand{\weg}[1]{}
\newtheorem{Cor}{Corollary}

\newcommand{\be}{\begin{equation}}
\newcommand{\ee}{\end{equation}}

\title[Projective equivalence of Randers metrics.]{On projective equivalence  and pointwise projective relation of Randers metrics.   } \date{}
\author{Vladimir S. Matveev}\thanks{Institute of Mathematics,  Friedrich-Schiller-Universit\"at Jena,  07737 Jena Germany\\  vladimir.matveev@uni-jena.de} 

\begin{document}


\begin{abstract} 
We correct a mistake in \cite{china} and prove the natural generalization of the projective Lichnerowicz-Obata conjecture for Randers metrics. 
\end{abstract} 

\maketitle

\section{Introduction}

\subsection{Definition and results} 
A Randers  metric  is a  Finsler metric of  the form 
\begin{equation} \label{a0}  \begin{array}{rl} F(x,\xi ) =& \sqrt{g(x)_{ij} \xi^i \xi^j}  + \omega(x)_{i}\xi^i \\ =&  \sqrt{g(\xi, \xi)}    +  \omega(\xi),\end{array}  
\end{equation} 
where $g= g_{ij}$ is a Riemannian metric  and $\omega= \omega_i$ is an $1$-form.  
 Here and everywhere in the paper we assume summation with respect to repeating indices. The assumption that $F$ given by \eqref{a0} is indeed a Finsler metric is equivalent to the condition that the $g$-norm of $\omega$ is less than one. 
Within the whole paper  we assume that all objects we consider are at least $C^2-$smooth, 
 
 By a {\it forward   geodesic}  of a Finsler metric  $F$ we understand a  regular curve $x:I\to M$ such for  any sufficiently close points 
 $a, b\in I$, $a\le b$  the restriction 
of the curve $x$ to the interval $[a,b]\subseteq  I$  
  is an extremal of the forward-length functional 
\begin{equation} 
\label{1bis}
 L^+_F(c):=  \int_{a}^b F(c(t), \dot c(t)) dt\end{equation} 
  in the set of all smooth curves $c:[a,b]\to M$ connecting $x(a)$ and $x(b)$.
 By a {\it backward   geodesic}  of a Finsler metric  $F$ we understand a  regular curve $x:I\to M$ such for  any sufficiently close points 
 $a, b\in I$, $a\le b$  the restriction 
of the curve $x$ to the interval $[a,b]\subseteq  I$  
  is an extremal of the backward-length functional 
\begin{equation} 
\label{2}
 L^-_F(c):=  \int_{a}^b F(c(t), -\dot c(t)) dt\end{equation} 
  in the set of all smooth curves $c:[a,b]\to M$ connecting $x(a)$ and $x(b)$.
 
 Note that these definitions do not assume any preferred parameter on the  geodesics: if $x(\tau)$ is a say forward geodesic, and $\tau(t)$ is a (orientation-preserving) 
  reparameterisation with $\dot \tau := \tfrac{d \tau }{dt}>0$, then $x(t):= x(\tau(t))$ is  also a forward  geodesic.  As the examples show, the  condition that the reparameterisation is orientation-preserving, i.e., the condition that 
  $\dot \tau := \tfrac{d \tau }{dt}>0$, is important though.

     We will study the question when two Randers metrics $F$ and $\bar F$ 
are {\it projectively equivalent}, that is, when every forward-geodesic of $F$   is a forward-geodesic of $\bar F$.  Within our paper we  will always  assume that the dimension is at least two, since in dimension one all metrics are projectively equivalent.

\begin{remark} \label{1} Comaparing \eqref{1bis} and \eqref{2} we see  that 
for every forward geodesic   $x(t)$, $t\in [-1,1]$,  
 the `reverse'  curve $\tilde x(t):= x(-t)$ is a backward-geodesic, and vice versa. Thus, if two metrics $F$ and $\bar F$ have the same   forward-geodesics, they automatically have the same backward-geodesics.  
 
Moreover, the forward-geodesics of any metric $F$ are the backward-geodesic of the metric $\tilde F$ given by  $\tilde F(x, \xi)=F(x, -\xi)$ and vice versa.  For the  Randers metrics, the transformation 
$F\mapsto \tilde F$ reads $\sqrt{ g(\xi,  \xi )} + \omega( \xi) \mapsto \sqrt{ g(\xi,  \xi) } - \omega( \xi)$.
     \end{remark}

There exists the following  `trivial'  examples of projectively equivalent Randers metrics: 

Every  Finsler metric $F$ is projectively equivalent to the Finsler  metric $\const \cdot F$ for any constant $\const>0$. Indeed,  forward and backward  geodesics are extremals of \eqref{1bis}, \eqref{2}. Now, replacing $F$ by $\const \cdot F$    multiplies  the length-functionals $L^{pm}$ 
of all curves by $\const$  so  extremals remain extremals. If the Finsler 
metric  $F$ is Randers, the operation $F\mapsto \const \cdot F$ multiplies the metric $g$ by $\const^2 $ and the form $\omega$ by $\const$. 

Any  Finsler metric  $F$ is   projectively equivalent to the Finsler metric $F+ \sigma$, where $\sigma$ is an arbitrary  closed one-form such that $F(x,\xi) +\sigma(\xi) > 0 $ for all tangent vectors $\xi\ne 0$. Indeed, 
a geodesic  connecting two points $x,y$  is an  extremal of the length functionals  
$L^\pm(c)$ given by \eqref{1bis}, \eqref{2},
over all regular curves  $c:[a,b]\to M$  with $c(a) = x$ and $c(b) = y$. Adding $\sigma$ to the Finsler metric $F$ changes the (forward and backward) length of all such curves  in one  homotopy class 
by adding a  constant  to it  (the constant may depend on the homotopy class), which does not affect the property of a curve to be an extremal. 
If the Finsler 
metric  $F$ is Randers, the operation $F\mapsto F + \sigma$ does not change the metric  $g$   and 
adds the form $\sigma$   to the form $\omega$.

One can combine two examples above as follows:
 any Finsler metric $F$ is projectively equivalent to the Finsler metric  
$\const \cdot F + \sigma$, where $\sigma$ is a closed 1-form and $\const \in \mathbb{R}_{>0}$ such that $\const\cdot  F(x,\xi) + \sigma(\xi)>0$ for all tangent vectors $\xi\ne 0$.
 
 Now, if the forms $\omega$ and $\bar \omega$  are closed, the Randers metrics 
 $F(x,\xi)=  \sqrt{ g( \xi,  \xi)}   +  \omega(\xi)$ and  $\bar F(x,\xi)=  \sqrt{ \bar g( \xi,  \xi)}   +  \bar\omega(\xi)$ on one manifold $M$  are projectively equivalent
if and only if the Riemannian metrics $g$ and $\bar g$ are projectively equivalent.  Indeed, as we explained above the geodesics of $F$ are   geodesics of the Riemannian metric $g$, the geodesics of $\sqrt{ \bar g(\xi, \xi)}  + \bar \omega(\xi)$ are the geodesics of the Riemannian metric $\bar g$, so that projective  equivalence of $ F$  and $\bar F$ is equivalent to 
the projective  equivalence of $ { g}$  and ${ \bar g}.$

Note that there are a lot of examples of projectively equivalent Riemannian metrics; the first examples were known already to Lagrange \cite{Lagrange}  and  local classification  of projectively equivalent metrics was  known already to Levi-Civita \cite{Levi-Civita}.

One of the goals of this note  is to show that the `trivial'  examples above give us all possibilities  of projectively equivalent Randers metrics.

\begin{Th}\label{thm2}  Let the Finsler metrics $\sqrt{g_{ij}\dot x^i \dot x^j}  + \omega_{i}\dot x^i$ and $\sqrt{\bar g_{ij}\dot x^i \dot x^j}  + \bar \omega_{i}\dot x^i$ on a connected manifold  be  projectively equivalent. Suppose at least one of the forms $\omega$ and $\bar \omega$  is not closed. Then,  for a certain $\const \in \mathbb{R}_{>0}$ we have 
$g = \const^2 \cdot   \bar g$ and the form  $\omega - \const\cdot  \bar \omega$ is closed. 
\end{Th}

 Let us first remark  that Theorem \ref{thm2} follows from \cite[Theorem 1.1]{Burns}.  
 The paper \cite{Burns} deals with the magnetic systems and study the question when magnetic 
 geodesics of one  magnetic system  $(g, \Omega)$ are reparameterized magnetic geodesics of another magentic system  $(\bar g, \bar \Omega)$,  see \cite{Burns} for definitions. 
In particular, it was    proved that if for positive numbers (energy levels) $E, \bar E\in \mathbb{R}$, 
every magnetic  geodesic with  energy $E$ of  one  magnetic system, is, after a proper reparameterisation, 
a magnetic geodesic with  energy $\bar E$ of another magnetic system, 
then $\bar g = \const\cdot  g$ and $\bar \Omega = \overline{\const} \cdot \Omega$ (the second constant $\overline{\const}$ depends on  $\const, E,  \bar E$), or $\Omega= \bar \Omega= 0$ and the metrics $g$ and $\bar g$ are projectively equivalent.   Now,  it is well known that forward geodesics of the  Randers metric \eqref{a0} are, after an appropriate orientation preserving reparameterisation,  magnetic  geodesics with energy $E=1$  of the  magnetic system $(g, \Omega = d\omega)$. In view of this, 
 Theorem \ref{thm2} is actually a corollary of \cite[Theorem 1.1]{Burns}.

 Theorem \ref{thm2} is visually  very close to \cite[Theorem 2.4]{china}: the only essential difference is that in the present paper we speak about \emph{projective equivalence}, and the condition discussed in 
  \cite{china}   is  when two metrics are \emph{pointwise projectively related} (see section \ref{mistake} for definition). In the Riemannian case (or, more generally, if the case when  Finsler metrics are reversible),  
  these two conditions, projective equivalence and pointwise projective relation,  coincide.  
    For general Finsler metrics, and in particular for Randers metrics,  projective equivalence and pointwise projective relation are different conditions  and    \cite[Theorem 2.4]{china}  is wrong:   we give a counterexample  in section \ref{mistake}. We also discuss how one can modify Theorem \ref{thm2} and  \cite[Theorem 2.4]{china} such that they  become correct, see Corollaries \ref{thm3} and \ref{cor4} from section \ref{mistake}.
 
 Besides, the proof of \cite[Theorem 2.4]{china}, even if one replaces pointwise projective relation projective equivalence, seems to have a certain  mathematical gap, namely  an important delicate   and nontrivial step  was not done (at least we did not find the place where it was discussed); out proor of Theorem \ref{thm2} closes this gap.  We comment on this in section \ref{mistake}.

As we mentioned above, we do not pretend that Theorem \ref{thm2} is new since it is a direct corollary of \cite[Theorem 1.1]{Burns}, though it seems to be unknown within Finsler geometers.  
New results of the papers  are related to projective transformations. By  {\it  projective transformation}  of $(M,F)$ we understand a  diffeomorphism  $\phi$ 
such that pullback of $F$ is projectively equivalent to $F$.  By  {\it homothety} of $(M,F)$ we understand a  diffeomorphism $\phi$  such that t pullback of $F$ is proportional  to $F$.  Homotheties evidently send forward  geodesics to forward  geodesics and are therefore  projective transformations.  As a direct application of  Theorem \ref{thm2}, we obtain 

\begin{Cor} \label{cor1}  If the form $\omega$ is not closed, every  projective transformation of \eqref{a0}  on a connected manifold $M$    is a   homothety of the Riemannian metric $g$. In particular, if $M$ is closed, every projective transformation of \eqref{a0} is an isometry  of the Riemannian metric $g$.
\end{Cor}

Many papers study  the question when  a  Randers metric is projectively flat, i.e., when  its forward  geodesics  are straight lines  in a certain coordinate system. Combining Theorem \ref{thm2}  with the classical Beltrami Theorem (see e.g. \cite{Beltrami,short,schur}), we obtain the following wellknown statement

\begin{Cor}[Folklore]  \label{cor2}  The metric \eqref{a0}  is projectively flat if and only if  $g$ has constant sectional curvature and $\omega$ is  closed. 
\end{Cor} 

In the case the manifold $M$ is closed (= compact and without boundary),  more can be said. 
We denote by $Proj(M,F)$  the group of the projective transformations of the Finsler manifold 
$(M,F)$ and by $Proj_0(M,F)$ its connected component containing the identity.

\begin{Cor} \label{cor3} Let $(M,F)$ be a closed connected Finsler manifold with $F$ given by \eqref{a0}.   Then, at least   one of the following possibilities holds:  
  
  \begin{enumerate} \item There exists a 
   closed form $\hat \lambda$  such that  $Proj_0(M,F)$ consists of isometries of the 
    the Finsler metric  $F(x,\xi)=  \sqrt{g(\xi,\xi) } + \omega(\xi) - \hat \lambda(\xi)  $, or 
  \item  the form  $\omega$ is  closed  and  $g$ has constant positive  sectional curvature. \end{enumerate} 
\end{Cor} 

\subsection{Motivation}   One of our    motivation  was to correct  mistakes in the paper \cite{china}:  to construct  counterexamples and to formulate the correct statement.  
 A part of this goal is to give a simple selfcontained proof of Theorem \ref{thm2} which  does not require Finsler machinery and therefore could be interesting for a bigger group of mathematicians.

 Actually,  
 a lot of papers discuss  metrics that are pointwise  projectively  related, and many of them have the same mistake as \cite{china}: in the proofs, the authors actually use that metrics are projectively equivalent, but formulate the results assuming the metrics are pointwise projectively related. 
 If for every forward geodesic $x:[-1, 1]\to M$ the reverse curve $x(-t)$ is also  a forward geodesic,  (for example, when the metrics are reversible),    then the results  remain correct; but in the general case many papers on pointwise   projectively  related metrics are wrong and in many cases in the papers it is not even mentioned whether the authors speak about all  metrics or restrict themself to the case when the metrics are reversible.   Since the Randers metrics are  nonreversible, this typical mistake is clearly seen in the  case of Randers metrics   and we have chosen \cite{china} to demonstrate it.

 Besides, we think that    Corollary 3 is deserved to be publishes, since  it is a natural generalisation of the classical projective Lichnerowicz-Obata conjectures for Randers metrics, see   \cite{nagano,Yamauchi1,hasegawa,solodovnikov1} where the Riemannian version of the conjecture was formulated and proved under certain additional geometric assumptions, and  \cite{obata,CMH,archive} where the Riemannian version of the conjecture was proved in the full generality.

Additional motivation to study projective equivalence and projective transformations came from mathematical relativity and lorentz differential geometry: it was observed that  the light-line geodesics of  a stationary, standard spacetime can be described with the help of Randers metrics on a manifold of dimension one less. This observation is  called the Stationary-Randers-Correspondence and  it is nowadays 
 a hot topic in  the Lorentz  differential geometry since one can effectively  apply it, see  for example  \cite{1,2,2a,how}.  The projective transformations of the Randers metrics correspond then  to the conformal  transformations of the initial Lorentz metric preserving the integral curves of the 
  Killing vector field, 
 so one can directly apply   our results.

\subsection{ Projective equivalence versus pointwise projective relation.} \label{mistake} 

By  \cite{china}, 
two Finsler metrics $F$ and $\bar F$  one one manifold $M$ 
are {\it pointwise  projectively   related}, if they have
the same geodesics as point sets.  The difference between this definition and our  definition of projective equivalence is that in our definition we also  require  that the orientation  of the forward 
geodesics is the same in both metrics.   In particular, the metrics $F$ and $\tilde F$  from  Remark \ref{1},    
such that every forward geodesic of the first is a backward geodesic of the second,    are pointwise  projectively related   according to the definition from  \cite{china}. This allows one to construct immediately a counterexample to \cite[Theorem 2.4]{china}.

\begin{example} Take any Riemannian metric $g$ and any  form $\omega$   such that it is not closed. Next, 
consider the Randers  metrics 
$F(x,\xi) = \sqrt{g(\xi, \xi)}  + \omega(\xi)$  and 
$\bar F(x,\xi) = \sqrt{g(\xi, \xi)}  - \omega(\xi)$ are pointwise 
 projectively  related,  
since for every forward  geodesic of $F$   is a backward geodesic of $\bar F$ and every backward   geodesic of $F$   is a forward  geodesic of $\bar F$.  Since the Riemannian parts   of these Finsler metrics coincide,   \cite[Theorem 2.4]{china} claims that the form $\omega - ( - \omega)= 2 \omega $ is closed 
 which is not the case. 
\end{example}

The following two corollaries   are  an attempt to correct the statement of   \cite[Theorem 2.4]{china}

\begin{Cor} \label{thm3} 
Suppose two Randers metrics  $ F(x,\xi) = \sqrt{g(\xi, \xi)}  +\omega(\xi)$ and  $  \bar F(x,\xi) = \sqrt{\bar g(\xi, \xi)}  +\bar \omega(\xi)$ on one connected manifold $M$ 
are pointwise projectively related. Assume in addition that  the set $M^0$  of the points of $M$  such that 
 the differential $d\omega$ is not zero is connected.
  Then, there exists   a  positive $\const\in \mathbb{R}$  such that  at least one of following statements holds at all points of the manifold: \begin{enumerate} \item \label{(1)}
    $\bar g   = \const^2 \cdot g$   and $\bar \omega - \const \cdot \omega$ is a closed form, or \item 
    \label{(2)} $\bar g   = \const^2\cdot g$   and $\bar \omega + \const \cdot \omega$ is a closed form.\end{enumerate}
\end{Cor}

It is important though that the  set $M^0$ is connected: Indeed, 
as the following example shows, the cases \eqref{(1)}, \eqref{(2)}  of Theorem \ref{thm3} could hold simultaneously  in different regions  of one manifold. 

\weg{\begin{example} \label{ex2}  Take $M= \mathbb{R}^2$ with the flat metric $g= d^2 {x} +  d^2{y}$. Next, consider the sequence of the points $P_1,...,P_k,...\in M$ such that it converges to $(0,0)$. Next, take the sequence of  radii $r_1,..., r_k,... >0 $ such that 
 the  closed  balls  $B_{r_k}(P_k)$ around the points  $P_1,...P_k,...$ with radii $r_1,..., r_k,... $  are disjunkt.  Consider a smooth form $\omega$ such that it
  vanishes outsides of the balls and its differential $d\omega$  is  not zero at least at one point of each ball.  Next, define 
   the form 
  $\bar \omega$    as follows: at the points of the ball $B_{r_k}(P_k)$ set 
  
  $$\left\{\begin{array}{rc}\bar \omega = \omega &  \textrm{if  $k$ is odd } \\  
   \bar \omega = -\omega &  \textrm{if  $k$ is even . }\end{array}\right.$$
At the other points set $\bar \omega=0$. 
 
    The form $\bar \omega$  is evidently smooth; the metrics $F$ and $\bar F$ are pointwise projectively related since  at the points of the balls  $B_{r_k}(P_k)$ with odd   $k$ every forward 
    geodesic of $F$ is a forward geodesic of $\bar F$; at the points of the balls  $B_{r_k}(P_k)$ with even   $k$ every forward     geodesic of $F$ is a backward  geodesic of $\bar F$; at all other 
     points  the metrics coincide. 
     
   By construction, in no neighborhood of the  point $(0,0) = \lim_{k\to \infty } P_k$ one of the cases listed in Theorem \ref{thm3} holds.   
\end{example} }

\begin{example} \label{ex2}  Consider the ray $S:= \{(x,y)\in \mathbb{R}^2\mid  \ x=0 \textrm{ and } y\le 2\}$ .  Take $M= \mathbb{R}^2 \setminus S$ 
with the flat metric $g= d^2 {x} +  d^2{y}$. Consider  two balls $B_+$ and  $B_-$ around the points 
$(+1,0)$ and $(-1,0)$ of radius $\tfrac{1}{2} $. Next, take a smooth  form $\omega $   on $M$ such that $\omega$ vanishes on $M\setminus \left( B_+ \cup B_-\right)$,  such that there exists at least one point of every ball such that $d\omega\ne 0$, and such that the $g$-norm of $\omega$ is less than $1$ at every point. Next, define 
   the form 
  $\bar \omega$    as follows:  set 
  $$\left\{\begin{array}{rc}\bar \omega = -\omega &  \textrm{at  $p\in B_-$  } \\  
   \bar \omega = \omega &  \textrm{at $p\not\in B_-$. }\end{array}\right.$$

    The form $\bar \omega$  is evidently smooth; the metrics $F(x,\xi)= \sqrt{g(\xi,\xi)} + \omega(\xi)  $ and $\bar F(x,\xi)= \sqrt{g(\xi,\xi)} + \bar \omega(\xi)$ are pointwise projectively related though none of the  of the cases listed in Corollary  \ref{thm3} holds on the whole manifold.   
\end{example}

\begin{Cor} \label{cor4}  Suppose two Rander metrics $F(x,\xi)= \sqrt{g(\xi,\xi)}+ \omega(\xi)$  and  
$\bar F(x,\xi)= \sqrt{ \bar g(\xi,\xi)}+ \bar \omega(\xi)$ on a connected  manifold  are pointwise projectively related. Suppose the form $\omega$ is not closed. Then, there exists a positive 
constant $\const \in \mathbb{R}$ such that   $g = \const^2 \cdot \bar g$  and such that 
 for every point 
 $x\in M$  we have $d\omega = \const \cdot d\bar \omega$ or  $d\omega = -\const \cdot d\bar \omega.$
\end{Cor} 

Let us now explain, as announced in the introduction,  one more mathematical difficulty with the proof of  \cite[Theorem 2.4]{china}.  Authors proved (in our notation and assuming that they actually work with projectively equivalent metrics) \begin{itemize} \item 
 that  at the points where the metrics $g$ and $\bar g$ are proportional with possibly a nonconstant coefficient,  the coefficient of proportionality is actually a constant, 
 $g = \const^2  \cdot \bar g$,  and   the form $\omega- \const\cdot  \bar \omega$ is   closed and  \item  that
   at the points  such that  the forms $\omega$ and $\bar \omega$ are closed the metric $g$ and $\bar g$   are projectively equivalent. \end{itemize} These two observations  do not immediately imply that one of these two conditions holds  on the whole manifold or  even, if they work locally, at all  points of   a  sufficiently small neighborhood  of arbitrary  point. Indeed,  
  we could conceive of a two Randers  metrics  $ \sqrt{g(\xi, \xi)} + \omega(\xi)$  
  and  $ \sqrt{\bar g(\xi, \xi)} + \bar \omega(\xi)$    such that at certain points of the  manifold 
   the metrics  $g$ and $\bar g$ 
   are  projectively equivalent but nonproportional (the  set of such points  is open), 
   and at  certain points the metrics  $g$ and $\bar g$ 
   are   proportional (the set of such points is evidently close).

We overcome 
 this difficulty by using the (nontrivial) result of  \cite[Corollary 2]{dedicata}, which implies that if two projectively equivalent metrics are proportional on a certain open set then they are proportional on the whole manifold (assumed connected). It is not the only possibility to overcome this difficulty, but at the present point we do not know any of them which is completely trivial; since this difficulty is not addressed in the proof  of  \cite[Theorem 2.4]{china} we suppose that the authors simply overseen the difficulty; our paper closes this gap.

\section{Proofs.}

\subsection{ Proof of Theorem \ref{thm2}.}  Recall that every  forward geodesic $x(t)$    of a 
 metrics is  an   extremal  of  the lengh functional $L(c)  = \int_a^b F(c(t), \dot c(t)) dt$, 
 and, therefore, is  a  solutions of the Euler-Lagrange equation 
 
 \begin{equation}\label{a1}
 \frac{d}{dt} \frac{\partial F}{\partial \dot x} - \frac{\partial F}{\partial x}=0.\end{equation} 
 For the  Randers metric \eqref{a0}, the equation \eqref{a1} reads 
 
 \begin{equation}  \label{a2} \begin{array}{cccl}  g_{ip} \ddot x^p \  \left( \frac{1}{\sqrt{g_{k m} \dot x^k \dot x^m}}   \right)
& +& g_{ip} \dot x^p \  \frac{d}{dt} \left( \frac{1}{\sqrt{g_{k m} \dot x^k \dot x^m}}   \right) +  \left( \frac{1}{\sqrt{g_{k m} \dot x^k \dot x^m}}   \right) \frac{\partial g_{ip}}{\partial x^k}  \dot x^k \dot x^p & \\ & +& \frac{\partial \omega_i}{\partial x^k} \dot x^k - \frac{\partial \omega_k}{\partial x^i} \dot x^k    - 
  \frac{1}{2 \sqrt{g_{k m} \dot x^k \dot x^m}} \dot x^p \dot x^q \frac{\partial g_{pq} }{\partial x^i}&=0.\end{array}
 \end{equation} 
 
 Multiplying this equation by $g^{ij}$ (the inverse matrix to $g_{ij}$), we obtain

 \begin{equation}  \label{a2bis} \begin{array}{cccl}   \ddot x^j \  \left( \frac{1}{\sqrt{g_{k m} \dot x^k \dot x^m}}   \right)
& +&  \dot x^j \  \frac{d}{dt} \left( \frac{1}{\sqrt{g_{k m} \dot x^k \dot x^m}}   \right) 
+  g^{ij} \left( \frac{1}{\sqrt{g_{k m} \dot x^k \dot x^m}}   \right) \frac{\partial g_{ip}}{x^k}  \dot x^k \dot x^p & \\ 
& +& g^{ij} \left(\frac{\partial \omega_i}{\partial x^k}  - \frac{\partial \omega_k}{\partial x^i}\right) \dot x^k    - 
g^{ij}  \frac{1}{2 \sqrt{g_{k m} \dot x^k \dot x^m}} \dot x^p \dot x^q \frac{\partial g_{pq} }{\partial x^i}&=0.\end{array}
 \end{equation}

 It is easy to check  by  calculations and is evident geometrically that for every solution $x(\tau)$  of \eqref{a2} 
 and for every time-reparameterization $\tau(t)$ with $\dot \tau > 0$ the curve $x(\tau(t))$ is also a forward geodesic. Thus,  if $\sqrt{  g_{ij}\dot x^i \dot x^j}  +  \omega_{i}\dot x^i$ and $\sqrt{ \bar g_{ij}\dot x^i \dot x^j}  + \bar \omega_{i}\dot x^i$ are projectively  equivalent,  
      then every solution $x(t)$ of  \eqref{a2}  also  satisfies 
 
 \begin{equation}  \label{a3bis} \begin{array}{cccl}   \ddot x^j \  \left( \frac{1}{\sqrt{\bar g_{k m} \dot x^k \dot x^m}}   \right)
& +&  \dot x^j \  \frac{d}{dt} \left( \frac{1}{\sqrt{\bar g_{k m} \dot x^k \dot x^m}}   \right) 
+  \bar g^{ij} \left( \frac{1}{\sqrt{\bar g_{k m} \dot x^k \dot x^m}}   \right) \frac{\partial \bar g_{ip}}{\partial x^k}  \dot x^k \dot x^p & \\ 
& +& \bar g^{ij} \left(\frac{\partial \bar \omega_i}{\partial x^k}  - \frac{\partial \bar \omega_k}{\partial x^i}\right) \dot x^k    - 
\bar g^{ij}  \frac{1}{2 \sqrt{\bar g_{k m} \dot x^k \dot x^m}} \dot x^p \dot x^q \frac{\partial \bar g_{pq} }{\partial x^i}&=0,\end{array}
 \end{equation}
 where $\bar g^{ij}$  is the inverse matrix to $\bar g_{ij}$.
 We now multiply the equation \eqref{a3bis} by ${\sqrt{\bar g_{k m} \dot x^k \dot x^m}}$ and subtract the equation \eqref{a2bis} multiplied by ${\sqrt{ g_{k m} \dot x^k \dot x^m}}$  to obtain

 \begin{equation}    \begin{array}{cl}
   & \dot x^j \left( \sqrt{\bar g_{k m} \dot x^k \dot x^m}\ \frac{d}{dt} \left( \frac{1}{\sqrt{\bar g_{k m} \dot x^k \dot x^m}}   \right) -   \sqrt{ g_{k m} \dot x^k \dot x^m}\ \frac{d}{dt} \left( \frac{1}{\sqrt{ g_{k m} \dot x^k \dot x^m}}   \right) \right)\\+ &\dot x^k\left( \bar g^{ij} \sqrt{\bar g_{k m} \dot x^k \dot x^m}\left(\frac{\partial \bar \omega_i}{\partial x^k} 
 - \frac{\partial \bar \omega_k}{\partial x^i}\right)  
 -g^{ij}\sqrt{ g_{k m} \dot x^k \dot x^m}  \left(\frac{\partial  \omega_i}{\partial x^k} 
 - \frac{\partial  \omega_k}{\partial x^i}\right) \right)\\
+ & \dot x^k \dot x^p\left( \bar g^{ij}  \frac{\partial \bar g_{ip}}{\partial x^k} 
-g^{ij}  \frac{\partial g_{ip}}{\partial x^k} \right)   
   - \dot x^p \dot x^q \left( \bar g^{ij}  \frac{1}{2 } \frac{\partial \bar g_{pq} }{\partial x^i}   +  g^{ij}  \frac{1}{2 }\frac{\partial  g_{pq} }{\partial x^i}\right) =0 .
 \end{array} \label{rez} \end{equation} 

In the equation above, the following logic in rearranging the terms was used: consider a forward-geodesic $\tilde x$ such that $\tilde x(0)= x(0)$ and such that $\dot {\tilde x}(0)= - \dot x(0)$. Then, at $t=0$, the first lines of the equation \eqref{rez} and its analog for $\tilde x$ are proportional to $\dot x(0)= -\dot{ \tilde x}(0)$.  The second line of the equation \eqref{rez} is minus  its analog for $\tilde x$. The remaining third  line  of  the equation \eqref{rez} coincides  with   its   analog  for $\tilde x$.

Then, subtracting the equation \eqref{rez}  from its analog for $\tilde x$ at $t=0$   we obtain (at $t =0$)
$$
    -\dot x^j f(x(0),\dot x(0)) + \dot x^k\left( \bar g^{ij} \sqrt{\bar g_{k m} \dot x^k \dot x^m}\left(\frac{\partial \bar \omega_i}{\partial x^k} 
 - \frac{\partial \bar \omega_k}{\partial x^i}\right)  
 -g^{ij} \sqrt{ g_{k m} \dot x^k \dot x^m}\left(\frac{\partial  \omega_i}{\partial x^k} 
 - \frac{\partial  \omega_k}{\partial x^i}\right) \right) =0 , 
$$
where
 $$\begin{array}{rc}f(x(0),\dot x(0)) = & \sqrt{\bar g_{k m} \dot x^k \dot x^m}\ 
\left(\frac{d}{dt}_{|t=0} \left( \frac{1}{\sqrt{\bar g(x)_{k m} \dot x^k \dot x^m}}\right) -  
 \frac{d}{dt}_{|t=0} \left( \frac{1}{\sqrt{\bar g(\tilde x)_{k m} \dot{ \tilde x}^k \dot{\tilde  x}^m}}\right)
\right)\\-&\sqrt{ g_{k m} \dot x^k \dot x^m}\ 
\left(\frac{d}{dt}_{|t=0} \left( \frac{1}{\sqrt{ g(x)_{k m} \dot x^k \dot x^m}}\right) -  
 \frac{d}{dt}_{|t=0} \left( \frac{1}{\sqrt{g(\tilde x)_{k m} \dot{ \tilde x}^k \dot{\tilde  x}^m}}\right)\right)
\end{array}$$

Since this  equation is fulfilled for every  geodesic, for every tangent vector $v $  we have

\begin{equation} \label{la1} 
   v^j  f(v)  =    \sqrt{\bar K(v)}  \bar g^{ij}\bar L_{ik} v^k 
 - \sqrt{K(v)} g^{ij}  L_{ik} v^k .
\end{equation}
 In the above equation, $$K(v)= g_{pq} v^p v^q, \   \bar K(v)= \bar g_{pq} v^p v^q, \  L_{ik} = d\omega = \frac{\partial  \omega_i}{\partial x^k}  - \frac{\partial  \omega_k}{\partial x^i},
\bar L_{ik} = d\bar \omega = \frac{\partial \bar  \omega_i}{\partial x^k}  - \frac{\partial  \bar \omega_k}{\partial x^i}. $$  
  Note that the matrices $L_{ij} $ and $\bar L_{ij}$ are scew-symmetric. 

Let us view this equation  as  a system of algebraic equations in a certain tangent space $\mathbb{R}^n$; we will  show that for every   symmetric positive definite matrices $g_{ij}$ and $\bar g_{ij}$ and for every   scew-symmetric matrices $L$, $\bar L$ there exist the following possibilities only:

\begin{itemize}  \item The matrix $g_{ij}$ is proportional to the matrix $\bar g_{ij}$ and  the matrix $L_{ij} $ is proportional to  the matrix $\bar L_{ij} $ with the same coefficient of proportionality, or  
\item The matrices   $L_{ij} $ and $\bar L_{ij} $ are zero.  
\end{itemize} 

Indeed, we multiply the equation  \eqref{la1} 
by $v^p g_{jp} $, and using  that $v^p g_{jp}\bar g^{ij}\bar L_{ik} v^k = v^iv^kL_{ik}=0$ because 
of scew-symmetry of $L$, we obtain 
$$ 
   K(v)   f(v)  =    - \sqrt{\bar K(v)} g_{jp} \bar g^{ij} \bar L_{ik} v^k v^p .
$$

Similarly,  multiplying  the equation  \eqref{la1} 
by $v^p \bar g_{jp} $,  we obtain 
$$   \bar K(v)   f(v)  =     \sqrt{ K(v)} v^p \bar g_{jp} g^{ij} L_{ik} v^k .
$$

Combining the last two equations we obtain 
\begin{equation}  \label{la3} 
\bar K(v)^3 \ (g_{jp} \bar g^{ij} \bar L_{ik} v^k v^p)^2 =  K(v)^3 \ (\bar g_{jp}  g^{ij} L_{ik} v^k v^p)^2.  
\end{equation} 
Since the algebraic expression $K(v)$ is irreducible (over $\mathbb{R}$), if 
$g_{jp} \bar g^{ij} \bar L_{ik} v^k v^p \not \equiv 0$, the equation 
 \eqref{la3} implies that $K(v) = \alpha  \bar K(v)$ implying $\bar g = \alpha^2 \cdot g $ for a certain $\alpha>0$. 
 In this case,   the equation \eqref{la1} implies $\bar L_{ij} = \alpha  L_{ij}$.  

Thus, at every point $p$ of our manifold, we have one of the two above possibilities. If at all points the second possibility takes places, the 1-forms $\omega$ and $\bar \omega$ are closed and the Riemannian metrics $g$ and $\bar g $ are  geodesically equivalent. If at least at one point the first possibility takes place, then in every point of  a small neighborhood of the point the first possibility takes place, so that  the differential of the forms $d\omega$ and $ d \bar \omega$ are  proportional at every point of a small neighborhood, 
 $d \omega =  \alpha(x) d\bar \omega$. The function $\alpha $ is then a constant implying the Finsler metrics are homothetic.   

Then, the metrics   $g$ and  $\bar g$ are projectively  equivalent: indeed, in the neighborhood consisting of   the points where  the 1-forms $\omega$ and $\bar \omega$ are closed  they are projectively equivalent as we explained above, and at the neighborhoods  where the differential 
 of at least one of the forms $\omega, \bar \omega$ 
  is not closed they are even proportional. Now, by \cite[Corollary 2]{dedicata}, if the metrics are projectively equivalent everywhere and nonproportional in some  neighborhood, their restrictions  to every neighborhood are  nonproportional.  Thus, in this case $d\omega=d\bar \omega=  0$ at every point of the manifold. Theorem \ref{thm2}  is proved. 

\subsection{Proof  of Corollaries \ref{cor1},\ref{cor2},\ref{cor3}. }  Corollary \ref{cor1} follows immediately from 
Theorem \ref{thm2}: if   $\phi:M\to M$  
is a  local  projective transformation of the metric \eqref{a0} on a connected manifold $M$, 
 then,  by the definition of projective transformations, 
   the pullback $\phi^*(F)$  is projectively equivalent to $F$. Then, by Theorem \ref{thm2},  if the form 
$\omega$ is not closed, we have $\phi^*(g) = \const\cdot g$ implying $\phi$ is a homothety for $g$. Corollary \ref{cor1} is proved.

In order to prove Corollary \ref{cor2}, we will use 
 that the straight lines are geodesics of the standard flat Riemannian  metric which we denote by $g^{flat}$. Then, a projectively flat metric is projectively equivalent  to  the Randers metric 
 $F^{flat}(x,\xi):= \sqrt{g^{flat}(\xi, \xi)}.  $ Then, by Theorem \ref{thm2},  the 
  form $\omega$ is  closed,  and the Riemannian metric $g$ is projectively flat. By  the classical Beltrami Theorem (see e.g. \cite{Beltrami,short,schur}), the metric $g$ has constant sectional curvature.  Corollary 2 is proved. 
  
  Let us now  prove Corollary \ref{cor3}.  We assume that $(M,F)$ is a closed connected  Finsler manifold with $F$ given by \eqref{a0} and denote by $Proj(M,F)$ its group of projective transformations and by $Proj_0(M,F)$ the connected component of this group containing the identity. 
   
  Suppose first the form  $\omega$ is  closed. Then,    every 
   projective transformation of the Randers metric \eqref{a0}  is a projective transformation of  the Riemannian metric $g$.    Then, if $Proj_0(M,F)$ contains not only isometries,  the metric $g$ has constant positive sectional curvature by the Riemannian 
   projective Obata conjecture (proven in  \cite[Corollary 1]{archive},   \cite[Theorem 1]{CMH} and \cite[Theorem 1]{obata}) and we are done.

  Assume now the form $\omega$ is not closed. Then, by Theorem \eqref{thm2}, 
  every element of $Proj_0(M,F)$ is an isometry of $g$ which in particular implies  that the group $Proj_0(M,F)$  is compact.   We consider an invariant measure $d\mu= d\mu(\phi)$ on $Proj_0(M,F)$ normalized such that \begin{equation} \int_{\phi\in Proj_0(M,F)}d\mu(\phi)= 1.\label{cl} \end{equation} 
  
   Consider  the 1-form $\hat \omega$ given by the formula 
  $$
  \hat\omega(\xi)= \int_{\phi\in Proj_0(M,F)}\phi^*\omega(\xi) d\mu(\phi).  
  $$
Here  $\phi^*\omega$ denotes the pullback of the form $\omega$ 
 with respect to the isometry $\phi\in Proj_0(M,F)$.

   Let us show that the form 
  $ \hat\omega - \omega$ is closed. In view of \eqref{cl},  we have 
  \begin{equation}
  \omega(\xi)- \hat \omega(\xi)  = \int_{\phi\in Proj_0(M,F)}(\omega - \phi^*\omega)(\xi)  d\mu(\phi).\label{cl2}\end{equation}

   Since for every $\phi\in Proj_0(M,F)$  the form $\lambda^\phi:= \omega - \phi^*\omega$ is closed by Theorem \ref{thm2}, the  form $\hat \lambda := \omega - \hat \omega$  is also closed. Indeed, let $\lambda^\phi_i$ be  the components of the form $ \lambda^\phi:= \omega - \phi^*\omega$  in a local 
   coordinate    system $(x^1,...,x^n)$.  Then, the components  $\hat \lambda_i$ 
   of $\omega - \hat \omega $ are given by the formula 
   $$
   \hat\lambda_i=  \int_{\phi\in Proj_0(M,F)}\lambda^\phi_i  d\mu(\phi).
   $$
  Differentiating this equation w.r.t. $x^j$, we obtain 
  $$ \tfrac{\partial}{\partial x^j} \hat\lambda_i=  \int_{\phi\in Proj_0(M,F)} \tfrac{\partial}{\partial x^j}  \lambda^\phi_i  d\mu(\phi).$$
  Now, since  for every $\phi\in Proj_0(M,F)$, 
   the form $ \lambda^\phi $ is closed, we have  $\tfrac{\partial}{\partial x^j}  \lambda^\phi_i = \tfrac{\partial}{\partial x^i}  \lambda^\phi _j$ implying $\tfrac{\partial}{\partial x^j} \hat\lambda_i= \tfrac{\partial}{\partial x^i} \hat\lambda_j$ implying $\hat \lambda $ is a closed form. 
   
   By construction,  the form $\hat \lambda$  satisfies the property that  $\omega - \hat \lambda$  is invariant with respect to $Proj_0(M,F)$, so the group $Proj_0(M,F)$ consists of isometries of the Finsler metric $\hat F(x, \xi)= \sqrt{g(\xi, \xi)} + \omega(\xi) - \hat \lambda(\xi)$.   Corollary \ref{cor3}  is proved. 

\subsection{Proof of Corollaries \ref{thm3}, \ref{cor4}.}

We will prove the Corollaries  \ref{thm3}, \ref{cor4} simultaneously. Within this section we assume that 
$F(x,\xi)= \sqrt{g(\xi, \xi)} + \omega(\xi)$ and 
$F(x,\xi)= \sqrt{\bar g(\xi, \xi)} + \bar \omega(\xi)$ are pointwise projectively related metrics on  a connected manifold $M$ and we denote by $M^0$ (resp.  $\bar M^0$)
 the set of the points of $M$ where the  differential of $\omega$ (resp. of $\bar \omega$) does not vanish.  $M^0$ and   $\bar M^0$ are  obviously open. 
 
Let us first observe that $M^0 =  \bar M^0$. In order to do it, consider  $p\in M^0$ and a vector $\xi \in T_pM$ such that $d\omega(\xi, \cdot)$ (viewed as an 1-form) does not vanish. Then, the forward and backward geodesic segments  $x(t)$ and $  \tilde x(t)$ such that $x(0)= \tilde x(0)= p$ and $\dot x(0)= \dot{\tilde x}(0)= \xi$ are two  different (even after a reparameterisation) curves. The curves are tangent at the point $x=0$. 

Indeed, in the proof of Theorem \ref{thm2} we found an ODE  \eqref{a2bis} for forward geodesics. Analogically, one can find an ODE for backward geodesics: it is  similar to \eqref{a2bis}, the only difference is that the term 
  $ g^{ij} \left(\frac{\partial\omega_i}{\partial x^k}  - \frac{\partial  \omega_k}{\partial x^i}\right) \dot x^k $ comes with the minus  sign.  We see that the difference between 
  $\ddot x(0)$ and $\ddot {\tilde x}(0)$ is  (in a local coordinate system) 
  a  vector which is not proportional to $\dot x(0)= \dot{\tilde x}(0)= \xi$. Then,  in a local coordinate system, 
  the  geodesic segments have different curvatures\footnote{In oder to define a curvature, we need to fix an euclidean structure in the neighborhood of $p$. The curvature depends on the euclidean structure, but the property of curvatures (considered as  vectors orthogonal to $\xi$) 
   to be different does not depend on the choice of the euclidean structure} 
     at the point $x(0)$. Then, the curves $x(t)$ and $\tilde x(\tau)$ do not have intersections for small $t\ne 0$, $\tau \ne 0$.  Moreover, would $d\omega(\xi, \cdot)= 0$, the second derivatives  $\ddot x(0)$ and $\ddot {\tilde x}(0)$ would coincide implying the curvatures of $x(t)$ and $\dot x(t)$ coincide at $t=0$.

     Now consider the geodesics of $\bar F$. There must be one (forward or backward) geodesic of $\bar F$ 
      that, after an orientation preserving  reparameterisation, 
        coincides   with $x(t)$, an one geodesic of $\bar F$ 
      that, after an orientation preserving  reparameterisation, 
        coincides   with $\tilde x(t)$. Since the curvatures of these two geodesics are different at the point $p$,  
   $d\bar \omega(\xi,\cdot)   \ne 0$ at the point $p$. Finally, an arbitrary point $p\in M^0$       also lies in $\bar M^0$, so that $ M^0\subseteq \bar M^0.$ Similarly one proves $ M^0\supseteq \bar M^0, $ so $M^0 $ and $\bar M^0$ coincide.

Let us now show that, on $M^0$, $d\omega =\const \cdot d\bar \omega$. In fact, we explain how one can modify the proof of Theorem \ref{thm2} to obtain this result.  As we explained above,  the forward and backward geodesic segments $  x(t)$ and $  \tilde x(t)$ such that $x(0)= \tilde x(0)$ and
 $\dot x(0)= \dot{\tilde x}(0)$ with $d\omega(\dot x(0), \cdot)\ne 0$ are geometrically different for small $t\ne 0$; then one of them is a forward geodesic of $\bar F$ and another is a backward geodesic of $\bar F$. We call $\xi\in T_pM, \ p \in M^0 $  a $p$-positive point, if the forward geodesic $x(t)$ with   $x(0)= p$ and
 $\dot x(0)= \xi$ is a forward geodesic of $\bar F$, and $p$-negative otherwise. 
 Next,  we consider two subsets of $T_pM$: 

$$S_+:= \{ \xi \in  T_pM\mid d\omega(\xi, \cdot )\ne 0;  \  d\bar \omega(\xi,\cdot)\ne 0; \ \textrm{$\xi $ is $p$-positive}\},$$
 $$
S_-:= \{ \xi \in  T_pM\mid d\omega(\xi,\cdot )\ne 0;  \  d\bar \omega(\xi, \cdot )\ne 0; \ \textrm{$\xi $ is $p$-negative}\}.$$ 
The closure 
 of one of  these two subsets contains a nonempty  open subset  $U\subset T_pM$; let us first assume that 
  the closure of $S_+$  contains a nonempty  open subset of $T_pM$.  Then, as in the proof of Theorem \ref{thm2}, we obtain that \eqref{la3} is valid for every  $v\in S_+$ (here we use a trivial fact that if $v \in S_+$ then $-v\in S_+$). Since \eqref{la3} is an algebraic condition, it must be then valid for all  $v\in T_pM$. Then, arguing as in 
  the proof of   Theorem \ref{thm2}, we conclude that   at the point $p$ we have 
 $\bar g = \alpha^2 \cdot g $  and $d\bar \omega   = \alpha\cdot   d\omega $ for some positive  $\alpha$.  
 Now, if   the closure of $S_-$  contains a nonempty  open subset of $T_pM$, we obtain  
 $\bar g = \alpha^2 \cdot g $  and $d\bar \omega   =  \alpha\cdot   d\omega $ for some negative  $\alpha$. 
 At $M^0$,  $\alpha$ must be  a smooth nonvanishing  function as the coefficient of the proportionality of two nonvanishing tensors; then  the sign of $\alpha$ is the  same of at all points of  every connected component of $M^0$. Since the forms  $d\bar \omega$ and $     d\omega $ are closed, the coefficient of the proportionality, i.e., $\alpha$ is  a constant. 
 
 This implies that  the Riemannaian  metrics $g$ and $\bar g$ are projectively equivalent, since at the points of $M^0$ they are even  proportional, and at the  inner points of $M\setminus M^0$ the geodesics of $g$ are geodesics of $F$ and geodesics of $\bar g$ are geodesics of $\bar F$.  By  \cite[Corollary 2]{dedicata}, the Riemannian  metrics $g$ and $\bar g$ are proportional at all points of the manifold so that $\bar g = \const^2 \cdot g$ on the whole manifold, as we claim in Corollaries \ref{thm3}, \ref{cor4}. 
  Comparing this with the condition 
 $\bar g = \alpha^2 \cdot g $  and $d\bar \omega   = \alpha\cdot   d\omega $  proved above for all points of each  connected component of $M^0$,   we obtain that at every connected component 
 of $M^0$  we have $d\bar \omega   = + \const\cdot   d\omega $  or $d\bar \omega   = - \const\cdot   d\omega $  (the sign can be different in  different  connected components of $M^0$  as Example \ref{ex2}  shows) as we claim in Corollaries  \ref{thm3}, \ref{cor4}. At the points of $M\setminus M^0$ we have $d\bar \omega   =   d\omega =0$  so both conditions $d\bar \omega   = + \const\cdot   d\omega $  and $d\bar \omega   = - \const\cdot   d\omega $ are fulfilled. Corollaries \ref{thm3}, \ref{cor4} are proved.

\end{document}